\magnification 1200
\pretolerance=500 \tolerance=1000  \brokenpenalty=5000
\mathcode`A="7041 \mathcode`B="7042 \mathcode`C="7043
\mathcode`D="7044 \mathcode`E="7045 \mathcode`F="7046
\mathcode`G="7047 \mathcode`H="7048 \mathcode`I="7049
\mathcode`J="704A \mathcode`K="704B \mathcode`L="704C
\mathcode`M="704D \mathcode`N="704E \mathcode`O="704F
\mathcode`P="7050 \mathcode`Q="7051 \mathcode`R="7052
\mathcode`S="7053 \mathcode`T="7054 \mathcode`U="7055
\mathcode`V="7056 \mathcode`W="7057 \mathcode`X="7058
\mathcode`Y="7059 \mathcode`Z="705A
\def\spacedmath#1{\def\packedmath##1${\bgroup\mathsurround =0pt##1\egroup$}
\mathsurround#1
\everymath={\packedmath}\everydisplay={\mathsurround=0pt}}
\def\nospacedmath{\mathsurround=0pt
\everymath={}\everydisplay={} } \spacedmath{2pt}
\def\qfl#1{\buildrel {#1}\over {\longrightarrow}}
\def\phfl#1#2{\normalbaselines{\baselineskip=0pt
\lineskip=10truept\lineskiplimit=0truept}\nospacedmath\smash
{\mathop{\hbox
to 10truemm{\rightarrowfill}}
\limits^{\scriptstyle#1}_{\scriptstyle#2}}}
\def\hfl#1#2{\normalbaselines{\baselineskip=0truept
\lineskip=10truept\lineskiplimit=1truept}\nospacedmath\smash{\mathop{\hbox
to12truemm{\rightarrowfill}}\limits^{\scriptstyle#1}_{\scriptstyle#2}}}
\def\gfl#1#2{\normalbaselines{\baselineskip=0truept
\lineskip=10truept\lineskiplimit=1truept}\nospacedmath\smash{\mathop{\hbox to
12truemm{\leftarrowfill}}\limits^{\scriptstyle#1}_{\scriptstyle#2}}}
\def\diagramme#1{\def\normalbaselines{\baselineskip=0truept 
\lineskip=10truept\lineskiplimit=1truept}   \matrix{#1}}
\def\vfl#1#2{\llap{$\scriptstyle#1$}\left\downarrow\vbox to
6truemm{}\right.\rlap{$\scriptstyle#2$}}

\def\mono{\lhook\joinrel\mathrel{\longrightarrow}}
\def\iso{\vbox{\hbox to .8cm{\hfill{$\scriptstyle\sim$}\hfill}
\nointerlineskip\hbox to .8cm{{\hfill$\longrightarrow $\hfill}} }} 

\def\sdir_#1^#2{\mathrel{\mathop{\kern0pt\oplus}\limits_{#1}^{#2}}}
\def\pprod_#1^#2{\raise
2pt \hbox{$\mathrel{\scriptstyle\mathop{\kern0pt\prod}\limits_{#1}^{#2}}$}}

      \font\eightbboard=msbm8
      \font\sevenbboard=msbm7 at 6pt
\font\fivebboard=msbm7 at 5pt
        
  \newfam\bboardfam
\font\eightrm=cmr8         \font\eighti=cmmi8
\font\eightsy=cmsy8        \font\eightbf=cmbx8
\font\eighttt=cmtt8        \font\eightit=cmti8
\font\eightsl=cmsl8        \font\sixrm=cmr6
\font\sixi=cmmi6           \font\sixsy=cmsy6
\font\sixbf=cmbx6\catcode`\@=11
\def\eightpoint{%
  \textfont0=\eightrm \scriptfont0=\sixrm \scriptscriptfont0=\fiverm
  \def\rm{\fam\z@\eightrm}%
  \textfont1=\eighti  \scriptfont1=\sixi  \scriptscriptfont1=\fivei
    \textfont2=\eightsy \scriptfont2=\sixsy \scriptscriptfont2=\fivesy
  \textfont\itfam=\eightit
  \def\it{\fam\itfam\eightit}%
  \textfont\slfam=\eightsl
  \def\sl{\fam\slfam\eightsl}%
  \textfont\bffam=\eightbf \scriptfont\bffam=\sixbf
  \scriptscriptfont\bffam=\fivebf
   \textfont\bboardfam=\eightbboard \scriptfont\bboardfam=
\sevenbboard\scriptscriptfont\bboardfam=\fivebboard
  \scriptscriptfont\bboardfam=\sevenbboard
  \def\Bbb{\fam\bboardfam}%
  \def\bf{\fam\bffam\eightbf}%

  \textfont\ttfam=\eighttt
  \def\tt{\fam\ttfam\eighttt}%
  \abovedisplayskip=9pt plus 3pt minus 9pt
  \belowdisplayskip=\abovedisplayskip
  \abovedisplayshortskip=0pt plus 3pt
  \belowdisplayshortskip=3pt plus 3pt 
  \smallskipamount=2pt plus 1pt minus 1pt
  \medskipamount=4pt plus 2pt minus 1pt
  \bigskipamount=9pt plus 3pt minus 3pt
  \normalbaselineskip=9pt
  \setbox\strutbox=\hbox{\vrule height7pt depth2pt width0pt}%
  \normalbaselines\rm}\catcode`\@=12

\newcount\noteno
\noteno=0
\def\up#1{\raise 1ex\hbox{\sevenrm#1}}
\def\note#1{\global\advance\noteno by1
\footnote{\parindent0.4cm
\up{\number\noteno}\
}{\vtop{\eightpoint\baselineskip12pt\hsize15.5truecm\noindent
#1}\vskip-5pt}\parindent 0cm}
\font\san=cmssdc10
\font\psan=cmssdc10 at 7pt
\def\pw{\kern1pt\hbox{\psan \char3}\kern1pt}
\def\ext{\hbox{\san \char3}}
\def\sym{\hbox{\san \char83}}

\def\pc#1{\tenrm#1\sevenrm}
\def\tx{\kern-1.5pt -}
\def\cqfd{\kern 2truemm\unskip\penalty 500\vrule height 4pt depth 0pt width
4pt\medbreak} 
\def\virg{\raise
.4ex\hbox{,}}
\def\decale#1{\smallbreak\hskip 28pt\llap{#1}\kern 5pt}
\def\no{n\up{o}\kern 2pt}
\def\ind{\par\hskip 1truecm\relax}
\def\indp{\par\hskip 0.5truecm\relax}
\def\moins{\mathrel{\hbox{\vrule height 3pt depth -2pt width 6pt}}}
\def\pmoins{\mathrel{\hbox{\kern1pt\vrule height 2pt depth -1pt
width 4pt\kern1pt}}}

\def\rond{\kern 1pt{\scriptstyle\circ}\kern 1pt}
\def\Hom{\mathop{\rm Hom}\nolimits}

\def\Ker{\mathop{\rm Ker}\nolimits}

\def\Pic{\mathop{\rm Pic}\nolimits}
\def\Tor{\mathop{\rm Tor}\nolimits}
\def\P{{\Bbb P}}
\def\C{{\Bbb C}}
\def\Z{{\Bbb Z}}
\def\pt{\scriptscriptstyle\bullet}
\def\cli{\mathop{\rm Cliff}\nolimits}
\frenchspacing
\input amssym.def
\input amssym
\input xy
\xyoption{all}
\vsize = 25truecm
\hsize = 16truecm
\voffset = -0.3truecm
\parindent=0cm
\baselineskip15pt
\overfullrule=0pt
\headline{\ifnum\pageno=1 \vtop{S\'eminaire BOURBAKI\hfill
Novembre 2003\par
 56\`eme ann\'ee, 2003-2004, \no 924
}
        \else \ifnum\pageno<10 {\hfill{\rm 924-0\folio}\hfill}
        \else {\hfill{\rm 924-\folio}\hfill} \fi \fi}
\nopagenumbers
\null
 \vskip 1.8truecm
\centerline{\bf LA CONJECTURE DE
GREEN G\'EN\'ERIQUE} 

\par \centerline{\bf [d'apr\`es C. Voisin]} \par \medskip
\centerline{par {\bf Arnaud BEAUVILLE}}  
\vskip 1.3truecm
{\bf 1. \'Enonc\'e de la conjecture}
\smallskip
\ind La conjecture de Green est une vaste g\'en\'eralisation de deux
r\'esultats classiques de la th\'eorie des courbes alg\'ebriques.  Soit $C$
une courbe complexe\note{Les th\'eor\`emes 1 et 2 ci-dessous sont
vrais en toute caract\'eristique [S-D]. Ce n'est pas le cas de la
conjecture de Green d'apr\`es [S1].}  projective et lisse (connexe), de
genre
$g\ge 2$. Soit
$K_C$ le fibr\'e canonique (= fibr\'e cotangent) de $C$.
On associe \`a
$C$ son {\it anneau canonique}
$$R:=\ \sdir_{n\ge 0}^{}H^0(C,K_C^n)\ .$$ 
\ind Notons $S$ l'alg\`ebre sym\'etrique
$\sym^{\scriptscriptstyle\bullet}H^0(C,K_C)$; c'est un anneau de
polyn\^omes en $g$ ind\'etermin\'ees.
\smallskip 
{\pc TH\'EOR\`EME} 1 (M. Noether)$.-$ {\it 
L'homomorphisme naturel
$S\rightarrow R$ est surjectif, sauf si $C$ est
hyperelliptique}.\smallskip 
\ind Supposons d\'esormais que $C$ n'est pas hyperelliptique. \`A
l'homomorphisme $S\rightarrow R$ correspond un plongement de
$C$ dans l'espace projectif $\P^{g-1}:=\P(H^0(C,K_C)^*)$, dit  {\it
plongement canonique}, qui  joue un r\^ole fondamental dans l'\'etude
de la g\'eom\'etrie de $C$. 
L'\'etape suivante est d'essayer de comprendre les \'equations de $C$
dans $\P^{g-1}$,  c'est-\`a-dire les \'el\'ements de $S$ qui s'annulent sur
l'image de $C$; ils forment un id\'eal gradu\'e $I_C$ de $S$, qui est
le noyau  de l'homomorphisme $S\rightarrow R$.
\smallskip 
{\pc TH\'EOR\`EME} 2 (Petri)$.-$ {\it L'id\'eal gradu\'e $I_C$ est
engendr\'e par ses \'el\'ements de degr\'e $2$, sauf si
$C$ est trigonale\note{La courbe $C$ est
dite trigonale si elle admet un morphisme $C\rightarrow \P^1$
de degr\'e 3.} ou isomorphe \`a une courbe plane de
degr\'e
$5$}.\smallskip 
\ind Chacun de ces deux th\'eor\`emes d\'ecrit la structure du $S$\tx
module $R$ en termes de l'existence de certains syst\`emes lin\'eaires
sur la courbe $C$. Par exemple, le th\'eor\`eme de Petri se traduit
(sauf pour les exceptions mentionn\'ees dans l'\'enonc\'e) par une suite
exacte
$$ S(-2)^{b_1}\longrightarrow S
\longrightarrow R\rightarrow 0\ ,$$ o\`u l'on note comme d'habitude
$S(-p)$ le $S$\tx module $S$ muni de la graduation d\'ecal\'ee de $p$
crans vers la droite:  $S(-p)_i=S_{i-p}$. 
\ind Cette pr\'esentation est un (petit) bout de la {\it r\'esolution
minimale} $P_{\pt}$ du $S$\tx module $R$, dont on sait depuis
Hilbert qu'elle est de la forme
$$0\rightarrow P_{g-2}\longrightarrow P_{g-3}\longrightarrow
\cdots\longrightarrow P_0\longrightarrow R\rightarrow 0\ ,$$o\`u
chaque $P_i$ est une somme directe de modules $S(-p)$, et o\`u les
diff\'erentielles sont donn\'ees par des matrices \`a
coefficients homog\`enes de degr\'e $\ge 1$. La r\'esolution
minimale est unique \`a isomorphisme (non unique) pr\`es\note{Dans le
langage des faisceaux, il revient au m\^eme de consid\'erer une
r\'esolution ${\cal P}_{\pt}$ du ${\cal O}_{\P^{g-1}}$\tx module
${\cal O}_C$, o\`u chaque ${\cal P}_i$ est une somme directe de
faisceaux 
${\cal O}_{\P^{g-1}}(-p)$, et o\`u les
 diff\'erentielles sont donn\'ees par des matrices \`a
coefficients homog\`enes de degr\'e $\ge 1$.}. 
\ind La dualit\'e de Serre entra\^{\i}ne que le complexe
$\Hom^{}_S(P_{\pt}, S(-g-1))$, d\'ecal\'e de $(g-2)$ crans vers la
gauche, d\'efinit encore une r\'esolution minimale de $R$, donc est
isomorphe \`a
$P_{\pt}$. Supposons $C$ non hyperelliptique; on a alors $P_0=S$
(Th. 1), donc  
$P_{g-2}=S(-g-1)$, et on s'aper\c{c}oit qu'il reste tr\`es peu de degr\'es
possibles pour les termes $P_i$ interm\'ediaires. De mani\`ere pr\'ecise, un
argument \'el\'ementaire montre qu'il existe un entier $c\ge 1$ tel que
$P_{\pt}$ soit de la forme:
$$\nospacedmath\displaylines{  0\rightarrow
S(-g-1)\longrightarrow  
  S(-g+1)^{b_1}\longrightarrow  S(-g+2)^{b_2}\longrightarrow 
\ \cdots\ \longrightarrow  S(-g+c-1)^{b_{c-1}}\longrightarrow  \cr
  S(-g+c+1)^{b'_c}\oplus
S(-g+c)^{b''_c} \longrightarrow \quad  \cdots\quad  
\longrightarrow     S(-c-2)^{b'_c}\oplus S(-c-1)^{b''_c}
\longrightarrow \cr  S(-c)^{b_{c-1}} 
\longrightarrow   \quad  \cdots\quad \longrightarrow S(-3)^{b_2}
\longrightarrow    S(-2)^{b_1}
\longrightarrow    S \ .\quad 
}$$
\ind La structure de la r\'esolution minimale est donc
essentiellement\note{Les entiers $b_i\ (i\le c-1)$ ainsi que $b''_c$
sont d\'etermin\'es par
$c$ et $g$, mais pas  $b'_i$ ni $b''_i$ pour $i>c$: le premier cas
  o\`u l'on trouve deux valeurs distinctes est $g=7, c=3$ [S1].}
d\'etermin\'ee par l'entier
$c$. 
\ind L'autre volet des th\'eor\`emes 1 et 2
 porte sur  la pr\'esence de syst\`emes lin\'eaires sp\'eciaux sur $C$. 
Si $L$ est un fibr\'e en droites sur $C$, de degr\'e $d$, on note
$h^i(L)$ la dimension de $H^i(C,L)\ \,(i=0,1)$, et l'on pose
$\cli(L):=g+1-(h^0(L)+h^1(L))=$ $d-2h^0(L)+2$; cet invariant v\'erifie
la relation agr\'eable
$\cli(L)=\cli(K_C\otimes L^{-1} )$. On d\'efinit alors l'{\it indice de
Clifford} $\cli(C)$ de
$C$ comme le minimum des entiers $\cli(L)$ pour tous les fibr\'es en
droites $L$ sur $C$ avec $h^0(L)\ge 2$ et $0\le d\le g-1$.
Un th\'eor\`eme classique de Clifford affirme que cet indice est toujours
positif, et qu'il est nul si et seulement si $C$ est hyperelliptique; de
plus les courbes d'indice 1 sont exactement celles qui apparaissent
dans le th\'eor\`eme 2. Les th\'eor\`emes 1 et  2 admettent donc la
reformulation suivante: 
$$\cli(C)\ge 1\  \Longleftrightarrow \  c\ge 1\quad , \quad 
\cli(C)\ge 2\ 
\Longleftrightarrow \  c\ge 2\ ,$$ce qui conduit  naturellement \`a la
\smallskip 
{\pc CONJECTURE DE} {\pc GREEN} [G]$.-$ $c=\cli(C)$.
\vskip1truecm
{\bf 2. R\'esultats}\vglue2pt
\ind Dans l'appendice de [G], Green et Lazarsfeld prouvent l'in\'egalit\'e
$c\le\cli(C)$, \`a l'aide des propri\'et\'es de la cohomologie
de Koszul \'etablies par Green dans le m\^eme article (voir \S 3). Il  s'agit
donc de d\'emontrer l'in\'egalit\'e oppos\'ee, c'est-\`a-dire, vu ce qui 
pr\'ec\`ede, que la composante
$(P_{p})_{p+1}$ de degr\'e $p+1$ de
$P_{p}$, avec $p=g-1-\cli(C)$, est nulle. 
\ind Cette conjecture remarquable a  vite attir\'e l'attention des
g\'eom\`etres alg\'ebristes.
   Dans [S1] Schreyer la v\'erifie  pour $g\le 8$; il
observe aussi qu'elle est fausse en carac\-t\'e\-ristique $2$, d\'ej\`a pour les
courbes g\'en\'erales de genre 7. 
 Le ``cas suivant" de la conjecture, $\cli(C)\ge
3\ \Longleftrightarrow
\  c\ge 3$, a \'et\'e d\'emontr\'e (ind\'ependamment) par Voisin [V1]
(pour $g\ge 11$), puis  Schreyer [S3] en g\'en\'eral. Le cas des courbes
planes est trait\'e dans [Lo].
 Divers cas particuliers ou reformulations de la
conjecture apparaissent dans [E], [P-R], [S2], [T]...
\ind Claire Voisin vient
de r\'esoudre le cas particuli\`erement int\'eressant des
courbes {\it g\'en\'erales} de genre
$g$. Elles v\'erifient
$\cli(C)=[{g-1\over 2}]$ (voir l'Appendice), de sorte que l'\'enonc\'e
prend une forme
parti\-cu\-li\`erement simple ({\it  conjecture de Green g\'en\'erique}):
\smallskip 
{\pc TH\'EOR\`EME} 3 ([V2], [V3])$.-$  {\it Posons
$c=[{g-1\over 2}]$. Pour une courbe de genre $g$ g\'en\'erale, la
r\'esolution minimale de
$R$ est de la forme
\vskip-20pt
$$ 0\rightarrow S(-g-1)\rightarrow 
  S(-g+1)^{b_1}\rightarrow  \cdots \rightarrow 
S(-c-2)^{b_c}\rightarrow    S(-c)^{b_{c-1}}
\rightarrow   \cdots \rightarrow   S(-2)^{b_1}
\rightarrow   S 
$$\vskip-5pt si $g$ est impair,  et\vskip-20pt
$$\nospacedmath\displaylines{\quad  0\rightarrow
S(-g-1)\rightarrow 
  S(-g+1)^{b_1}\rightarrow  \cdots \rightarrow 
S(-c-3)^{b_c-1} \rightarrow  \hfill \cr\hfill
S(-c-2)^{b_c/2}\oplus S(-c-1)^{b_c/2}
\rightarrow  S(-c)^{b_c-1}
\rightarrow   \cdots \rightarrow   S(-2)^{b_1}
\rightarrow   S \quad 
}$$\vskip-10pt si $g$ est pair}.
\smallskip \ind En fait la m\'ethode de d\'emonstration donne un
r\'esultat plus fort. Pour des entiers $g$ et $p$ fix\'es, consid\'erons
l'ensemble des courbes de genre $g$ $p$\tx {\it gonales}, c'est-\`a-dire
admettant un morphisme de degr\'e $p$ sur $\P^1$. Elles sont
param\'etr\'ees par un sch\'ema irr\'eductible, le sch\'ema de Hurwitz.
Nous verrons au \S 6 qu'une variante de la d\'emon\-stration du
th\'eor\`eme 3 dans le cas  $g$ pair entra\^{\i}ne la conjecture de Green pour
les courbes   $p$\tx gonales assez  g\'en\'erales, pour $p\ge {g\over
3}+1$. Or il se trouve que  M. Teixidor a obtenu (par une m\'ethode tr\`es
diff\'erente) le  r\'esultat correspondant pour 
$p\le  {g\over 3}+2$ [T]. Ainsi:
\smallskip 
{\pc TH\'EOR\`EME} 4 ([V2], [T])$.-$  {\it Une
courbe   $p$\tx gonale  g\'en\'erale v\'erifie la conjecture de
Green}.\smallskip  
\ind Plus pr\'ecis\'ement, on a  $c=\cli(C)=p-2$ pour $p\le [{g+3\over
2}]$. 
L'int\'er\^et de cet \'enonc\'e vient de ce que pour presque toutes\note{Au
moins conjecturalement -- voir l'Appendice pour une formulation pr\'ecise.}
les courbes
$C$,  l'indice de Clifford est \'egal \`a $\gamma - 2$, o\`u $\gamma$ (la
``gonalit\'e") est le plus petit entier  tel que $C$ soit $\gamma$\tx
gonale.
\ind Signalons que le Th. 3 pour $g$ impair a la cons\'equence suivante, qui
avait \'et\'e observ\'ee par Hirschowitz et Ramanan avant la d\'emonstration
de [V3], et qui apporte un peu plus d'eau au moulin de la conjecture
de Green: 
\smallskip 
{\pc COROLLAIRE} [H-R]$.-$ {\it Supposons $g=2k-1$. Dans l'espace des
modules des courbes de genre $g$, le lieu des courbes qui n'ont pas la
r\'esolution minimale g\'en\'erique co\"\i ncide avec celui des courbes
$k$\tx gonales}.

\vskip1truecm
{\bf 3. Cohomologie de Koszul}\vglue3pt \ind 
 Consid\'erons plus g\'en\'eralement une vari\'et\'e projective $X$, munie
d'un faisceau ample $L$. Notons   \vskip-7pt
$$V=H^0(X,L)\qquad S=\sym^{\pt}V\qquad 
R=\ \sdir_{n}^{}H^0(X,L^n)\ ;$$on s'int\'eresse \`a la r\'esolution
gradu\'ee libre minimale $P_{\pt}$ du
$S$\tx module gradu\'e  $R$. 
Consid\'erons $\C$ comme un $S$\tx module
via l'homomorphisme d'augmentation $S\rightarrow \C$.  Le
$S$\tx module gradu\'e  
$\Tor_i^S(\C,R)$ se calcule en substituant \`a $R$ la r\'esolution
$P_{\pt}$; comme celle-ci est minimale, le complexe
$\C\otimes_S P_{\pt}$ est \`a diff\'erentielle nulle, et l'on trouve donc
des isomorphismes de $S$\tx modules gradu\'es
$\Tor_i^S(\C,R)\iso \C\otimes_S  P_i $. Mais on peut aussi calculer ce
module en utilisant une r\'esolution libre gradu\'ee de $\C$. Il en existe
une bien connue, le complexe de Koszul
$$0\rightarrow  \ext^nV \otimes_{\C}  S(-n)\longrightarrow
\ldots\longrightarrow   
\ext^2V \otimes_{\C} S(-2)  \longrightarrow    V\otimes_{\C} S(-1)
\longrightarrow S$$ (avec $n=\dim V$). La diff\'erentielle
$ d_p: 
\ext^pV\otimes_{\C}S(-p)\longrightarrow  
\ext^{p-1}V\otimes_{\C} S(-p+1) $ applique
$(v_1\wedge\ldots\wedge v_p)\otimes P\ $ sur  $\displaystyle
\ \sum_i (-1)^{i+1}
 (v_1\wedge\ldots\wedge v_{i-1}\wedge
v_{i+1}\wedge\ldots\wedge v_p)\otimes P.v_i$ .

\ind Ainsi la composante de degr\'e $p+q$ du $S$\tx module gradu\'e
$\Tor_p^S(\C,R)$ s'identifie \`a l'espace d'homologie ${\cal
K}_{p,q}(X,L)$  du complexe
$$\ext^{p+1}V\otimes R_{q-1}\ \hfl{d_{p+1}}{}\ \ext^pV\otimes
R_q\ \hfl{d_{p}}{}\ \ext^{p-1}V\otimes R_{q+1}\ .$$ 
\ind Les espaces ${\cal K}_{p,q}(X,L)$ (``cohomologie de Koszul")
poss\`edent un grand nombre de pro\-pri\'et\'es int\'eressantes, \'etudi\'ees
notamment dans [G]. L'une d'elles sera fondamentale pour ce qui suit:
supposons pour simplifier
$L$ tr\`es ample, de sorte que $X$ est plong\'ee dans un espace projectif
de fa\c{c}on que $L={\cal O}_X(1)$.
Soit $Y$ une
section hyperplane\note{Cela signifiera ici qu'aucune
composante de $X$ n'est contenue dans l'hyperplan $\ell =0$.} de
$X$,  d\'efinie par une \'equation
$\ell =0$ (avec
$\ell \in H^0(X,L)$). Consid\'erons les anneaux
$S_Y=\sym^{\pt}H^0(Y,L_{|Y})$ et
$R_Y=\ \sdir_{n}^{}H^0(Y,L_{|Y}^n)$. Faisons en outre l'hypoth\`ese
$H^1(X, L^i)=0$ pour tout $i\ge 0$; elle garantit que $S_Y$ s'identifie
\`a
$S/(\ell )$ et $R_Y$ \`a $R/(\ell )$. Si   $P_{\pt}$ est une
r\'esolution minimale du $S$\tx module $R$, alors $P_{\pt}/\ell
P_{\pt}$ est une r\'esolution minimale du $S_Y$\tx module $R_Y$.  On
en d\'eduit un {\it isomorphisme canonique} ${\cal K}_{p,q}(X,L)\iso
{\cal K}_{p,q}(Y,L_{|Y})$ (``th\'eor\`eme de Lefschetz").

\smallskip 
 \ind  Revenons \`a notre courbe $C$. D'apr\`es le d\'ebut du \S 2, 
 la conjecture de Green se
traduit  par l'annulation de ${\cal K}_{p,1}(C,K_C)$ pour
$p=g-1-\cli(C)$, ou encore par l'exactitude de la suite
$$\ext^{p+1}H^0(C,K_C)\ \hfl{d_{p+1}}{}\
\ext^pH^0(C,K_C)\otimes H^0(C,K_C)\ \hfl{d_{p}}{}\
\ext^{p-1}H^0(C,K_C)\otimes H^0(C,K_C^2)\ .$$
Pour une courbe $C$ g\'en\'erale de genre $g$, l'indice de
Clifford vaut $[{g-1\over 2}]$, et il s'agit donc de prouver  l'annulation
de ${\cal K}_{k,1}(C,K_C)$ avec $k=[{g\over 2}]$.  Il  
 suffit de l'obtenir pour {\it une} courbe de genre $g$; 
 C.~Voisin utilise des courbes situ\'ees sur des surfaces tr\`es
particuli\`eres, les surfaces K3.

\ind Rappelons que les surfaces K3 sont, par d\'efinition, les surfaces
(lisses, compactes) simplement connexes \`a fibr\'e canonique trivial.
Celles qui nous int\'eressent ici sont les surfaces K3 $X$ polaris\'ees de
genre $g$, c'est-\`a-dire munies d'un fibr\'e en droites tr\`es ample $L$ de
carr\'e $2g-2$; on supposera de plus que la classe de $L$ dans le groupe
de Picard $\Pic(X)$ n'est  divisible par aucun entier $\ge 2$. Les
sections globales de
$L$  d\'efinissent un plongement de
$X$ dans
$\P^g$, dans lequel  les  sections hyperplanes lisses de $X$ sont
des  courbes de genre $g$, plong\'ees dans $\P^{g-1}$ par
le plongement canonique. 
Pour chaque entier $g\ge 3$, les surfaces K3 polaris\'ees de genre
$g$ forment une
 famille irr\'eductible; une surface assez
g\'en\'erale\note{\kern-3pt C'est-\`a-dire situ\'ee en dehors d'une
r\'eunion d\'enombrable d'hypersurfaces dans l'espace des
param\`etres.} dans cette famille v\'erifie
$\Pic(X)=\Z[L]$. 
\ind Les sections hyperplanes d'une telle surface $X$
 ne sont pas g\'en\'eriques pour $g\ge 12$, mais elles tendent
\`a se comporter 
 comme la courbe g\'en\'erique, en particulier du point de vue
de la th\'eorie de Brill-Noether [L]: par exemple  leur
indice de Clifford est l'indice g\'en\'erique $[{g-1\over 2}]$. Il est donc
tout-\`a-fait naturel d'essayer de prouver l'annulation de ${\cal
K}_{k,1}(C,K_C)$, avec $k=[{g\over 2}]$, pour ces courbes.
D'apr\`es le ``th\'eor\`eme de Lefschetz" pour la cohomologie de Koszul,
elle est \'equivalente \`a l'annulation de
 ${\cal K}_{k,1}(X,L)$.   La courbe $C$ va
d\'esormais dispara\^{\i}tre au profit de la surface K3 $X$.

\vskip1truecm
{\bf 4. Le cas de genre pair: strat\'egie de la preuve}\smallskip 
\ind La premi\`ere id\'ee force de la d\'emonstration est l'interpr\'etation
de ${\cal K}_{p,1}(X,L)$ en termes du {\it sch\'ema de Hilbert}
 de
$X$. Si $X$ est une vari\'et\'e projective et $d$ un entier, le sch\'ema de
Hilbert $X_d$ (not\'e plut\^ot d'habitude $X^{[d]}$ ou ${\rm
Hilb}^d(X)$) param\`etre les sous-sch\'emas finis de longueur $d$ de
$X$. Rappelons qu'un tel sous-sch\'ema $Z$ consiste en la donn\'ee de
 points $x_1,\ldots ,x_m$ de $X$ et en chacun de ces points d'un
id\'eal ${\cal I}_{x_i}$ de l'anneau local ${\cal O}_{x_i}$, de fa\c{c}on
que $\sum_i\dim_{\bf C}({\cal O}_{x_i}/{\cal I}_{x_i})=d$. En
associant  \`a $Z$ l'ensemble $\{x_1,\ldots ,x_m\}$,
chaque $x_i$ \'etant compt\'e avec sa multiplicit\'e $\dim({\cal
O}_{x_i}/{\cal I}_{x_i})$, on obtient un morphisme birationnel
$\varepsilon$ de
$X_d$ sur la puissance sym\'etrique $d$\tx i\`eme $\sym^dX$. Lorsque
$X$ est une {\it surface}, $X_d$ est lisse et irr\'eductible, de sorte que
$\varepsilon$ fournit une r\'esolution des singularit\'es de $\sym^dX$.
Nous nous bornerons \`a ce cas dans la suite\note{La Prop. 1 ci-dessous
s'\'etend en toute dimension \`a condition de se limiter aux
sous-sch\'emas finis  {\it curvilignes}, c'est-\`a-dire contenus dans une
courbe lisse.}.
\ind Soit $I_d$ la sous-vari\'et\'e de $ X\times X_d$ form\'ee
des couples $(x,Z)$ tels que $x$ soit un point de
$Z$. C'est une vari\'et\'e normale, munie de projections:
$$\xymatrix{&I_d\ar[dl]_p\ar[dr]^q&\\
X && X_d\ \ ;\\}$$la
fibre de $q$ en un point $Z$ de $X_d$  
s'identifie au sous-sch\'ema $Z$ de $X$. 
\ind  On associe \`a tout fibr\'e en droites $L$ sur $X$ le fibr\'e vectoriel 
${\cal E}_L:=q_*(p^*L)$ sur $X_d$, de rang $d$; sa fibre en un point
$Z$ de $X_d$ s'identifie \`a $H^0(Z,L_{|Z})$. On pose $L_d:=\det {\cal
E}_L$. Une analyse pr\'ecise du fibr\'e en droites $q^*L_d$ conduit
alors  au r\'esultat suivant:
\smallskip 
 {\pc PROPOSITION} 1$.-$ {\it L'espace ${\cal
K}_{d-1,1}(X,L)$ s'identifie au conoyau de
l'homo\-mor\-phisme}
$ q^*:H^0(X_d,L_d)\rightarrow H^0(I_d,q^*L_d)$.
\smallskip 
\ind  La d\'emonstration sera esquiss\'ee au \S 5. Comme expliqu\'e au \S
3, le th\'eor\`eme 3 pour $g$ pair r\'esultera de:
\smallskip 
 {\pc PROPOSITION} 2$.-$ {\it Soit $X$ une surface} K3 {\it dont
le groupe de Picard est engendr\'e par un fibr\'e ample $L$, avec
$L^2=2g-2$ et $g=2d-2$. L'homomorphisme
$q^*:H^0(X_d,L_d)\rightarrow $ $H^0(I_d,q^*L_d)$ est
surjectif}.
\ind Comme le morphisme
$q$ est fini et plat, on dispose
d'un homomorphisme dans l'autre sens
$q_*:H^0(I_d,q^*L_d)\rightarrow H^0(X_d,L_d)$, qui v\'erifie
$q_*\rond q^*=d$. La surjectivit\'e de $q^*$ est donc
\'equivalente \`a l'injectivit\'e de $q_*$.  
\ind Le c\oe ur de la d\'emonstration consiste alors \`a construire une
vari\'et\'e $Z$, munie d'un morphisme $j:Z\rightarrow X_d$, telle
que  le carr\'e cart\'esien\vskip-10pt
$$\diagramme{\widetilde{Z}&\hfl{\tilde j}{}& I_d\cr
\vfl{q^{}_Z}{} & & \vfl{}{q}\cr
Z&\hfl{j}{} & X_d}$$\vskip-10pt
poss\`ede les deux propri\'et\'es suivantes:
\indp (i) l'homomorphisme $\tilde j^*: H^0(I_d,q^*L_d)\rightarrow
H^0(\widetilde{Z},\tilde j^*q^*L_d)$ est injectif;
\indp (ii)  l'homomorphisme
$(q^{}_Z)_* :H^0(\widetilde{Z},q_Z^*(j^*L_d))\rightarrow
H^0(Z,j^*L_d)$ est injectif.
\ind Au vu du diagramme commutatif\vskip-10pt
$$\diagramme{H^0(\widetilde{Z},\tilde j^*q^*L_d) &\gfl{\tilde
j^*}{}& H^0(I_d,q^*L_d) \cr
\vfl{(q^{}_Z)_*}{} & & \vfl{}{q_*}\cr
H^0(Z,j^*L_d)&\gfl{j^*}{} & H^0(X_d,L_d) }$$\vskip-10pt
on en d\'eduit aussit\^ot l'injectivit\'e de $q_*$, et donc la Proposition 2.
\smallskip 
\ind La construction de $Z$ repose sur l'existence  d'un fibr\'e
vectoriel $E$ de rang 2 remarquable sur $X$, introduit par Lazarsfeld
dans [L]. C'est l'unique fibr\'e
de rang 2 stable sur $X$ de d\'eterminant $L$ et seconde classe de Chern
$c_2(E)=d$; il v\'erifie $\dim H^0(X,E)=d+1$.
  En associant \`a une section $s$ de $E$ son sch\'ema des z\'eros
$Z(s)$, on d\'efinit un morphisme 
$\P(H^0(X,E))\rightarrow X_d$ (qui est d'ailleurs un plongement).
Notons
$W$ l'image r\'eciproque de  $\P(H^0(X,E))$ dans $I_d$.
Elle est form\'ee des couples $(s,x)$ dans $\P(H^0(X,E))\times X$ tels
que $s(x)=0$. Pour $(s,x)$ dans un ouvert convenable $W^{\rm o}$ de
$W$,  le sch\'ema r\'esiduel $Z(s)\moins x$ est bien d\'efini.
Consid\'erons
l'application rationnelle $j_0: X\times W^{\rm o}\dasharrow X_d$ qui
associe \`a
$(y,(s,x))$ le sch\'ema $(Z(s)\moins x)\cup y$. En \'eclatant dans
$X\times W$ le lieu des $(y,(s,x))$ tels que $y\in Z(s)\moins x$ et en
restreignant \`a un gros ouvert, on obtient le morphisme
$j:Z\rightarrow X_d$ cherch\'e.
\ind Le c\oe ur de la d\'emonstration consiste alors \`a v\'erifier les
propri\'et\'es (i) et (ii) ci-dessus. Cette v\'erification prend 30 pages tr\`es
denses de l'article [V2], qu'il n'est pas question de reproduire ici.
J'essaierai d'en indiquer quelques \'etapes au paragraphe suivant.
\vskip1truecm
{\bf 5. Le cas de genre pair: quelques d\'etails}\smallskip 
 a) {\it D\'emonstration de la Proposition} 1\vglue0pt
\ind  Suivant [V2], nous dirons qu'un ouvert $V^{\rm o}$
d'une vari\'et\'e normale $V$ est {\it gros} si le ferm\'e
compl\'ementaire est de codimension $\ge 2$. Si $L$ est un
fibr\'e sur $V$, l'application de restriction $H^0(V,L)\rightarrow
H^0(V^{\rm o},L)$ est alors un isomorphisme.

\ind La premi\`ere \'etape  est le calcul de
$H^0(X_d,L_d)$. Les homomorphismes de restriction
$H^0(X,L)\rightarrow H^0(Z,L_{|Z})$, pour
$Z\in X_d$, d\'efinissent une fl\`eche $H^0(X,L)\otimes
_{\C}{\cal O}_{X_d}\rightarrow {\cal E}_L$, d'o\`u en passant aux
$\ext^d$ un homomorphisme
$\ext^d H^0(X,L)\rightarrow H^0(X_d,L_d)$, qui est en fait un {\it
isomorphisme} : on le voit en rempla\c{c}ant $X_d$ par le gros ouvert
des sous-sch\'emas ayant au plus un point double, et en \'ecrivant ce
dernier comme quotient d'un gros ouvert de $X^d$ \'eclat\'e le long des
diagonales $x_i=x_j$. 
\ind  On va d\'esormais remplacer $X_d$ par le gros ouvert des
sous-sch\'emas {\it curvilignes}, c'est-\`a-dire contenus dans une
courbe lisse -- et $I_d$ par l'ouvert des couples $(x,Z)$ o\`u $Z$ est
curviligne.  Pour un tel couple le sch\'ema r\'esiduel $Z\moins x$ est
bien d\'efini; on
 dispose donc d'un morphisme $$\tau :I_d\rightarrow
X\times X_{d-1} \quad \hbox{d\'efini par}\quad  \tau (x,Z)=(x,Z\moins
x)\ .$$ C'est un isomorphisme sur l'ouvert $U$  de
$I_d$ form\'e des couples $(x,Z)$ pour lesquels $x$ est un point simple
de $Z$; il contracte le diviseur $D:=I_d\moins U$ sur la vari\'et\'e
d'incidence 
$I_{d-1}\i X\times X_{d-1}$.
\ind On d\'eduit facilement de la d\'efinition de $L_d$ un isomorphisme
$$q^*L_d\cong \tau ^*(L\boxtimes L_{d-1})(-D)\ ,\eqno(*)$$
d'o\`u une suite exacte:
$$0\rightarrow H^0(I_d,q^*L_d)\rightarrow H^0(X\times
X_{d-1},L\boxtimes L_{d-1})\rightarrow
H^0(I_{d-1}, (L\boxtimes L_{d-1})^{}_{\,|\,I_{d-1}})\ .$$  Notons
$\tau ': I_{d-1}\rightarrow X\times X_{d-2}$ l'application 
correspondant \`a $\tau$, et $D'$ le diviseur de $I_{d-1}$ contract\'e par
$\tau '$. Appliquant de nouveau $(*)$ on trouve un
isomorphisme $(L\boxtimes L_{d-1})^{}_{|I_{d-1}}\cong
\tau'^*(L^2\boxtimes L_{d-2})(-D')$, d'o\`u une injection de
$H^0(I_{d-1}, (L\boxtimes L_{d-1})^{}_{|I_{d-1}})$ dans $
H^0(X,L^2)\otimes H^0(X_{d-2}, L_{d-2})$. On a finalement une  suite
exacte:
$$0\rightarrow H^0(I_d,q^*L_d)\rightarrow H^0(X,L)\otimes
H^0(X_{d-1},L_{d-1})\qfl{\beta} H^0(X,L^2)\otimes
H^0(X_{d-2},L_{d-2})\ ,$$de sorte que le conoyau de
$q^*:H^0(X_d,L_d)\rightarrow H^0(I_d,q^*L_d)$ s'identifie \`a
l'homologie d'un complexe
$$H^0(X_d,L_d)\qfl{\alpha}H^0(X,L)\otimes
H^0(X_{d-1},L_{d-1})\qfl{\beta} H^0(X,L^2)\otimes
H^0(X_{d-2},L_{d-2}) \ ;$$on v\'erifie
 que ce complexe s'identifie via les isomorphismes 
$\ext^{p}H^0(X,L)\iso H^0(X_p,L_p)$ au complexe de Koszul
$$\ext^{d}H^0(X,L)\ \phfl{d_d}{}\ H^0(X,L)\otimes
\ext^{d-1}H^0(X,L)\ \phfl{d_{d-1}}{}\  H^0(X,L^2)\otimes
\ext^{d-2}H^0(X,L)\ ,$$d'o\`u la Proposition 1.

\medskip

b) {\it D\'emonstration de la Proposition$\,2:$ propri\'et\'e} (i)
\ind Dans la suite de ce paragraphe il est commode de poser
$k=d-1$ (de sorte qu'on a
$g=2k$). Notons  $\P$ le gros ouvert de $\P(H^0(X,E))$ form\'e des
sections dont le sch\'ema des z\'eros est curviligne. Reprenons le carr\'e
cart\'esien\vskip-10pt
$$\diagramme{W & \hfl{}{} & I_{k+1}&\cr
\vfl{\pi }{} &&\vfl{}{}&\cr
\P&\hfl{}{}&X_{k+1}&\kern-3pt ;}$$\vskip-10pt soit
$\psi:W\rightarrow X_{k}$ le morphisme $(\sigma
,x)\mapsto Z(\sigma )\moins x$. En explicitant la d\'efinition de $Z$ on
se ram\`ene facilement \`a prouver l'injectivit\'e de l'application
$$\psi^*:H^0(X_{k}, L_{k})\rightarrow H^0(W,\psi^*L_{k})\ .$$
\ind Le point cl\'e pour cela est la construction, \`a partir d'une \'etude
fine  du fibr\'e de Lazarsfeld $E$, d'un isomorphisme canonique
$\ext^{k}H^0(X,L)\iso \sym^kH^0(X,E)^*$. D'autre part on
montre que le fibr\'e
$\psi^*L_{k}$ est isomorphe \`a $\pi ^*{\cal
O}_{\P}(k)$, d'o\`u un homomorphisme
injectif
$\sym^{k}H^0(X,E)^*\mono H^0(W,\psi^*L_{k})$. On conclut
en  v\'erifiant que le diagramme\vskip-10pt
$$\diagramme{\ext^{k}H^0(X,L)&\hfl {}{}&\sym^kH^0(X,E)^*\cr
\vfl{}{} & &\vfl{}{} \cr
H^0(X_{k},
L_{k})&\hfl{\psi^*}{}& H^0(W,\psi^*L_{k})
}$$\vskip-10pt est commutatif \`a un scalaire pr\`es.\smallskip 
c) {\it
D\'emonstration de la propri\'et\'e} (ii)
\ind Notons $\widetilde{W}$ le produit fibr\'e
$W\times_{X_{k+1}}I_{k+1}$, de sorte qu'on a un carr\'e cart\'esien
\vskip-10pt
$$\diagramme{\widetilde{W} & \hfl{}{} & I_{k}&\cr
\vfl{q^{}_W}{} &&\vfl{}{q}&\cr
W& \hfl{\psi}{} &X_{k} &.}$$\vskip-10pt
 Apr\`es quelques 
p\'erip\'eties, on se ram\`ene \`a prouver la surjectivit\'e de
l'homomorphisme $q_W^*:H^0(W,\psi^*L_{k})\rightarrow
H^0(\widetilde{W}, q_W^*\psi^*L_{k})$. Rappelons qu'on a 
$\psi^*L_{k}\cong \pi ^*{\cal O}_{\P}(k)$. Notons $r$ 
l'application compos\'ee
$\widetilde{W}\rightarrow W\rightarrow \P$.
En fait Voisin prouve un
r\'esultat plus fort, \`a savoir:
\indp$\bullet$ {\it L'homomorphisme} $r^*:H^0(\P,{\cal
O}_{\P}(k))\rightarrow H^0(\widetilde{W},r^*{\cal
O}_{\P}(k)))$ {\it est surjectif}. 
\ind La d\'emonstration de ce
r\'esultat occupe 16 pages de [V2] et je ne peux faire mieux qu'y
renvoyer le lecteur. Disons simplement qu'on r\'ealise
$\widetilde{W}$ comme un sous-sch\'ema de $B_\Delta(S\times
S)\times \P$, o\`u $B_\Delta(S\times
S)$ est obtenu en \'eclatant $S\times S$ le long de la diagonale. La
surjectivit\'e cherch\'ee est  \'equivalente \`a l'annulation d'un $H^1$
convenable sur $B_\Delta(S\times S)\times \P$.
Des
calculs de cohomologie d\'elicats sur cette vari\'et\'e ram\`enent cette
annulation \`a des \'enonc\'es sur les sections globales du fibr\'e de Lazarsfeld.
 \font\itg=cmmib10
\vskip1truecm
{\bf 6. La conjecture de Green pour les courbes {\itg p}-gonales
g\'en\'erales} 
\smallskip \ind Soit toujours $X$ notre surface K3, munie d'un
fibr\'e en droites $L$ v\'erifiant
$L^2=2g-2$, $g=2d-2$ et $\Pic(X)=\Z\,[L]$. Comme promis,
nous allons voir que l'annulation de
${\cal K}_{d-1,1}(X,L)$ entra\^{\i}ne le Th\'eor\`eme 4.  Choisissons des
points g\'en\'eraux $x_1,\ldots ,x_\delta $ de $X$, avec $\delta \le
(d-1)/2$.  Comme $\dim H^0(X,E)=d+1$, il existe  deux sections
lin\'eairement ind\'ependantes $s,t $ de $E$ s'annulant en ces points.
Pour un choix g\'en\'erique  des $x_i$ et de $s,t$, la courbe
$C$  o\`u s'annule  la section
$s\wedge t$ de $\ext^2E=L$ est lisse sauf en
$x_1,\ldots ,x_\delta $, o\`u elle a des points doubles ordinaires. Soit
$n:N\rightarrow C$ sa normalisation.
\smallskip 
{\pc PROPOSITION} 3$.-$ {\it La courbe $N$ est $(d-\delta )$\tx
gonale, et v\'erifie} $\ c=\cli(N)=d-\delta -2$.
\ind La courbe $N$ est de genre $\gamma=2d-2-\delta $; les
in\'egalit\'es $0\le \delta \le (d-1)/2$ se traduisent par 
${\gamma\over
3}+1\le d-\delta \le {\gamma\over 2}+1$. Cela donne la conjecture
de Green pour les
 courbes $p$\tx gonales g\'en\'erales de genre $\gamma$, avec
${\gamma\over 3}+1\le p\le {\gamma\over 2}+1$, et donc, compte
tenu de [T], le th\'eor\`eme 4.
\smallskip 
{\it D\'emonstration de la Proposition} 3 : Les sections $s,t$
engendrent un sous-faisceau de rang 1 de $n^*(E_{|C})$; la
partie mobile du syst\`eme lin\'eaire correspondant est un pinceau de
degr\'e $d-\delta $ (le nombre de z\'eros de $s$ ou $t$ en dehors des
$x_i$). La courbe $N$ est donc $(d-\delta )$\tx gonale, et il suffit de
prouver qu'on a ${\cal K}_{p,1}(N,K_N)=0$ pour
$p=\gamma-1-(d-\delta -2)=d-1$. Or l'annulation de
${\cal K}_{d-1,1}(X,L)$ (Prop. 1 et 2) garantit celle de
${\cal K}_{d-1,1}(C,K_{C})=0$; il s'agit de comparer ${\cal
K}_{d-1,1}(C,K_{C})$ et ${\cal K}_{d-1,1}(N,K_{N})$. L'application
trace $n_*K_N\rightarrow  K_C$ fournit  des injections naturelles
$H^0(N,K_N)\mono H^0(C,K_C)$ et $H^0(N,K_N^2)\mono
H^0(C,K_C^2)$, d'o\`u un diagramme commutatif\vskip-4pt
$$\xymatrix@=20pt{\ \ext^dH^0(N,K_N) \
\ar[r]_-{d^{}_N}\ar[dd]^{j'}&
\ \ext^{d-1}H^0(N,K_N)\otimes H^0(N,K_N)\
\ar@{->}[r]\ar[dd]^j 
\ar@/_2pc/[l]^{r^{}_N}&
\ \ext^{d-2}H^0(N,K_N)\otimes H^0(N,K_N^2)\ \ar[dd]^{j''} \\ 
\\
\ \ext^dH^0(C,K_C)\  \ar[r]_-{d_C} &
\ \ext^{d-1}H^0(C,K_C)\otimes H^0(C,K_C)\  \ar[r]
\ar@/_2pc/[l]^{r^{}_C}&
\ \ext^{d-2}H^0(C,K_C)\otimes H^0(C,K_C^2)
}$$\smallskip 
Les diff\'erentielles $d^{}_N$ et $d^{}_C$ admettent des r\'etractions
canoniques $r^{}_N$ et $r^{}_C$, d\'efinies par $r_{\pt}(\tau \otimes
\omega)={1\over d}\,\omega\wedge\tau $, qui commutent aux
fl\`eches verticales; cela entra\^{\i}ne que {\it l'homomorphisme 
${\cal
K}_{d-1,1}(N,K_{N})\rightarrow {\cal K}_{d-1,1}(C,K_{C})$ induit
par $j$ est injectif}. En effet, si un \'el\'ement $v$ de
$\ext^{d-1}H^0(N,K_N)\otimes H^0(N,K_N)$ est tel que $j(v)$ est
un bord, on a
$$j(v)=d^{}_Cr^{}_Cj(v)=d^{}_Cj'r^{}_N(v)=jd^{}_Nr^{}_N(v)\
,$$d'o\`u, puisque $j$ est injectif, $ v=d^{}_Nr^{}_N(v)$. Ainsi ${\cal
K}_{d-1,1}(C,K_{C})$ est nul, d'o\`u la Proposition 3.
\vskip1truecm {\bf 7. Le cas de genre impair}
\vglue3pt 
\ind Ce qui pr\'ec\`ede repose de mani\`ere essentielle sur les propri\'et\'es
du fibr\'e de Lazarsfeld, qui n'existe qu'en genre pair. Pour traiter le
cas  $g$ impair, C. Voisin consid\`ere une surface K3 $X$ dont le
groupe de Picard est engendr\'e par un fibr\'e en droites tr\`es ample $L$,
de carr\'e $2g-2$,  et la classe d'une courbe rationnelle lisse $\Delta$
telle que $\deg (L_{|\Delta})=2$. Posons $L'=L(\Delta)$. On a
$L'^2=2g$, $\deg (L'_{|\Delta})=0$; le morphisme $X\rightarrow
\P^{g+1}$ associ\'e \`a $L'$ est un plongement en dehors de $\Delta$ et
contracte $\Delta$ sur un point. 
\ind Posons $g=2k +1$. La premi\`ere \'etape de la d\'emonstration est de
v\'erifier que la Proposition 2 s'\'etend \`a $(X,L')$,  donnant ${\cal
K}_{k+1,1}(X,L')=0$. La d\'emonstration de la propri\'et\'e (i) s'adapte
imm\'ediatement, celle de (ii) demande nettement plus de travail.
\smallskip 
\def\pt{\scriptscriptstyle\bullet}
\ind Il s'agit maintenant d'en d\'eduire  l'annulation de  ${\cal
K}_{k,1}(X,L)$. Il est  commode pour cela d'utiliser 
 la dualit\'e de Serre, qui fournit une dualit\'e canonique entre
${\cal K}_{p,1}(X,L)$ et ${\cal K}_{g-2-p,2}(X,L)$.  Ainsi
${\cal K}_{k-1,2}(X,L')$ est nul, et on veut en d\'eduire l'annulation de 
l'espace ${\cal K}_{k-1,2}(X,L)$. Rappelons que celui-ci
est l'homologie du complexe\note{Dans ce paragraphe, pour tout
faisceau $F$ sur $X$ on note simplement $H^0(F)$ l'espace
$H^0(X,F)$.}
$$\ext^{k}H^0(L) \otimes  H^0(L)\ \hfl{d_{L}}{}\
\ext^{k-1}H^0(L) \otimes  H^0(L^2)\ \hfl{d_{L}}{} \
\ext^{k-2}H^0(L) \otimes  H^0(L^3)\ .$$
\ind Au couple $(X,L')$ est associ\'e comme plus haut le fibr\'e de
Lazarsfeld $E$, de d\'eterminant $L'$; la preuve repose sur la
construction d'un homomorphisme
 $$\varphi
:\sym^kH^0(E)\rightarrow \ext^{k-1}H^0(L) \otimes
H^0(L(-\Delta))$$ qui s'inspire d'une construction analogue utilis\'ee par
Green et Lazarsfeld pour prouver l'in\'egalit\'e $c\le
\cli(C)$  ([G], Appendice). 
\'Etant donn\'e deux sections globales $v,w$ de $E$, on notera $v\pw w$
leur produit ext\'erieur dans
$\ext^2 E=L'$.  On choisit
une base
$(w_1,\ldots ,w_{k+1})$ de $H^0(E(-\Delta))$, et on pose\note{Le
chapeau sur un terme signifie comme d'habitude qu'on l'omet.}, pour
$v\in H^0(E)$,
$$\varphi (v^k)=\sum_{i<j}(-1)^{i+j}\,(v\pw 
w_1)\wedge\ldots\wedge\widehat{(v\pw
w_i)}\wedge\ldots\wedge\widehat{(v\pw w_j)}\wedge\ldots
\wedge(v\pw w_{k+1})\,\otimes \,(w_i\pw w_j)\ ;$$la
condition 
 $w_i\in H^0(E(-\Delta))$ entra\^{\i}ne bien $v\pw w_i\in H^0(L)$ et
$w_i\pw w_j\in H^0(L(-\Delta))$.
\ind Choisissons d'autre part une section $\sigma $ de $H^0(L')$
dont la restriction \`a $\Delta$ n'est pas nulle; elle fournit un scindage
de la suite exacte
$$0\rightarrow H^0(L)\rightarrow H^0(L')\rightarrow
H^0(L'_{\,|\Delta})\cong\C\rightarrow 0\ ,$$
d'o\`u une d\'ecomposition
$H^0(L')=H^0(L)\oplus
\C\,\sigma $. Consid\'erons le diagramme commutatif
$$\diagramme{\ext^{k-1}H^0(L) \otimes 
H^0(L(-\Delta)) & \hfl{\delta }{} & 
\ \ext^{k-2}H^0(L)\otimes H^0(L^2(-\Delta))\cr
\vfl{1\otimes \sigma }{}& & \vfl{}{1\otimes \sigma}\cr
\ext^{k-1}H^0(L) \otimes  H^0(L^2)\quad\ &\hfl{d_{L}}{} &
\ext^{k-2}H^0(L) \otimes  H^0(L^3)\ .\quad 
}$$
\ind L'annulation
de ${\cal K}_{k-1,2}(X,L)$ va r\'esulter des quatre points suivants:
\smallskip 
\indp(i) L'homomorphisme compos\'e $\delta \rond \varphi$ est nul.
\indp(ii) L'homomorphisme induit
$\varphi:\sym^kH^0(E)\rightarrow \Ker \delta $ est surjectif.
\indp(iii) ${\cal K}_{k-1,2}(X,L)$ est engendr\'e par les classes
d'\'el\'ements $(1\otimes \sigma )\cdot \alpha$ pour
$\alpha\in\Ker\delta $.
\indp(iv) Pour $t\in\sym^kH^0(E)$, la classe de $(1\otimes \sigma
)\cdot \varphi (t)$ dans ${\cal K}_{k-1,2}(X,L)$ est nulle.
\smallskip 
\ind Les assertions (i) et (iv) r\'esultent d'un calcul sans myst\`eres, bas\'e sur
l'identit\'e 
$$(v_1\pw v_2)\cdot (v_3\pw v_4)-(v_1\pw v_3)\cdot (v_2\pw
v_4)+(v_1\pw v_4)\cdot (v_2\pw v_3)=0\quad \hbox{dans
}H^0(L'^2)$$quels que soient $v_1,\ldots ,v_4$ dans $H^0(E)$. 
\ind Prouvons (iii). Soit $\beta\in\Ker d_L$. Puisque ${\cal
K}_{k-1,2}(X,L')=0$, il existe un \'el\'ement $\gamma$ de
$\ext^{k}H^0(L') \otimes  H^0(L')$ tel que
$\beta=d^{}_{L'}\gamma$. La d\'ecomposition
$H^0(L')=H^0(L)\oplus
\C\,\sigma $ permet d'\'ecrire\kern1.2cm
$\gamma=\gamma_1+\sigma \wedge
\gamma_2+\gamma_3\otimes \sigma +(\sigma \wedge
\gamma_4)\otimes \sigma\ ,\qquad \qquad  $avec\smallskip 

\centerline{$\gamma_1\in \ext^{k}H^0(L) \otimes  H^0(L)\ ,\
\gamma_2\in \ext^{k-1}H^0(L) \otimes  H^0(L)\ ,\
\gamma_3\in\ext^{k}H^0(L)\ ,\ \gamma_4\in\ext^{k-1}H^0(L)\ .
 $}\smallskip 
\ind L'\'el\'ement $\gamma_4$ s'identifie \`a
l'image de
$d^{}_{L'}\gamma$ dans $\ext^{k-1}H^0(L') \otimes 
H^0(L'^2_{|\Delta})$; comme $d^{}_{L'}\gamma=\beta$ appartient \`a
$\ext^{k-1}H^0(L) \otimes  H^0(L^2)$, on en d\'eduit
$\gamma_4=0$. Comme on peut modifier $\beta$ par un bord on
peut supposer $\gamma_1=0$. Enfin on a
$$\gamma_3\otimes \sigma =d^{}_{L'}(\sigma \wedge
\gamma_3)+\sigma \wedge d^{}_{L'}\gamma_3\ ,$$de sorte qu'en
modifiant $\gamma$ par un bord on peut supposer
$\gamma_3=0$. 
\ind On a alors $\gamma=\sigma \wedge \gamma_2$, et par suite
$\beta=d^{}_{L'}\gamma=\gamma_2\cdot (1\otimes \sigma)
-\sigma \wedge d^{}_{L}\gamma_2$. En utilisant de nouveau la
d\'ecomposition
$H^0(L')=H^0(L)\oplus
\C\,\sigma $,  le fait que $\beta$ appartient \`a
$\ext^{k-1}H^0(L) \otimes  H^0(L^2)$  implique d'une part que
$d^{}_{L}\gamma_2$ est nul, d'autre part que $\gamma_2$
appartient au sous-espace
$\ext^{k-1}H^0(L) \otimes  H^0(L(-\Delta))$. Par suite $\gamma_2$
appartient \`a $\Ker \delta $, et
 ${\cal K}_{k-1,2}(X,L)$ est engendr\'e par les
classes des \'el\'ements $\gamma_2\cdot (1\otimes \sigma)$ avec 
$\gamma_2\in \Ker \delta $.
\smallskip 
\ind Le gros du travail est la d\'emonstration de (ii), pour laquelle je
ne peux que renvoyer \`a [V3], p. 12-26. Disons simplement qu'on se
ram\`ene \`a un \'enonc\'e sur la cohomologie d'un \'eclatement convenable
de $\P(H^0(E))\times X$, \'enonc\'e dont la d\'emonstration demande
une ing\'eniosit\'e technique consid\'erable.
\vskip1truecm
{\bf Appendice: l'indice de Clifford}\smallskip 
\ind Nous utiliserons dans cet appendice une abr\'eviation tr\`es
classique: un  syst\`eme lin\'eaire\note{\'Etant donn\'e un diviseur
$D$, le syst\`eme lin\'eaire $|D|$ est l'ensemble des diviseurs effectifs
lin\'eairement \'equivalents \`a $D$; il s'identifie \`a l'espace projectif
$\P(H^0(C,{\cal O}_C(D))$.}
$|D|$ sur
$C$ de degr\'e
$d$ et de dimension projective $r$ est appel\'e un $g^r_d$. L'indice de
Clifford 
$\cli(C)$  est alors le minimum des entiers $d-2r$ sur l'ensemble des
$g^r_d$ avec $d\le g-1$ et $r\ge 1$. D'apr\`es le 
 th\'eor\`eme  de Clifford, on a $\cli(C)\ge 0$, et $\cli(C)=0$
 si et seulement si $C$ est hyperelliptique.  Cet invariant a \'et\'e
introduit par Martens dans [M], o\`u il montre entre autres que
 les courbes d'indice 1 sont exactement celles qui apparaissent
dans le th\'eor\`eme 2.
\ind Consid\'erons les courbes de genre $g$ fix\'e. Les courbes $d$\tx
gonales, c'est-\`a-dire admettant un $g^1_d$, ont un indice de Clifford
$\le d-2$; lorsqu'elles sont assez g\'en\'erales, leur indice de Clifford est
exactement $d-2$ [Ba]. Il s'ensuit que l'indice de Clifford  est
$[{g-1\over 2}]$ pour une courbe g\'en\'erale, et prend toutes les
valeurs entre $0$ et $[{g-1\over 2}]$. 
\ind  Les courbes dont l'indice de Clifford est fourni par
un $g^r_d$ avec $r>1$ (et pas par un syst\`eme lin\'eaire de
dimension plus petite) sont beaucoup plus rares. Pour $r=2$, ce sont
les courbes planes lisses de degr\'e $d$, qui sont de genre ${1\over
2}(d-1)(d-2)$. Pour $3\le r\le 9$, les auteurs de [ELMS] prouvent
que cela impose  $g=4r-2$, avec un
indice de Clifford 
$2r-3$ donn\'e par un fibr\'e en droites $L$ tel que $L^2\cong K_C$;
ils conjecturent bien naturellement  le m\^eme \'enonc\'e pour tout $r$
(et construisent, pour tout $r$, une courbe ayant les propri\'et\'es
indiqu\'ees). Si cette conjecture est correcte, et si $C$ est une courbe de
genre $g$ et d'indice de Clifford $c$, alors:
\indp a) $C$ est $(c+2)$\tx gonale, ou\indp b) 
$g={1\over
2}(c+2)(c+3)$, $C$ est une courbe plane lisse de degr\'e $c+4$, ou
\indp c)
$c$ est impair $\ge 3$, $g=2c+4$, et $C$ admet un fibr\'e en droites
$L$ tel que $L^2\cong K_C$ et $\cli(L)=c$. 
\ind Pour $g$ fix\'e il y a donc (moyennant la conjecture) au plus deux
valeurs de $c$ pour lesquelles il existe des courbes   d'indice $c$ qui
ne soient pas $(c+2)$\tx gonales.
\vskip1truecm{\eightpoint\baselineskip=12pt
\leftskip1cm\rightskip1cm\hskip0.8truecm Je remercie Olivier Debarre et
Claire Voisin pour leurs commentaires  pertinents sur une premi\`ere
version de ce texte.\par}

 \vskip 2truecm
\centerline{{\bf BIBLIOGRAPHIE}}
\vglue15pt\baselineskip12.8pt
\font\cc=cmcsc10
\def\num#1#2#3{\smallskip\item{\hbox to\parindent{\enskip [#1]
\hfill}}{\cc#2} : {\sl #3}}
\parindent=1.3cm 
\num{E}{L. Ein}{A remark on the syzygies of the generic canonical
curves}. J. Differential Geom. {\bf 26} (1987), 361--365.
\num{ELMS}{D. Eisenbud, H. Lange, G. Martens,
F.-O. Schreyer} {The Clifford dimen\-sion of a projective
curve}. Compositio Math. {\bf 72} (1989),  173--204.
\num{G}{M.  Green}{Koszul cohomology and the geometry of
projective varieties}. J. Differential Geom. {\bf 19} (1984),  125--171.
\num{H-R}{A. Hirschowitz, S. Ramanan}{New evidence for Green's
conjecture on sy\-zy\-gies of canonical curves}. Ann. Sci. \'Ecole Norm.
Sup. (4) {\bf 31} (1998),  145--152.
\num{L}{R. Lazarsfeld}
{Brill-Noether-Petri without degenerations}.
J. Differential Geom. {\bf 23} (1986),  299--307.
\num{Lo}{F. Loose}{On the graded Betti numbers of plane algebraic
curves}. Manuscripta Math. {\bf 64} (1989),  503--514.
\num{M}{H. Martens}{Varieties of special divisors on a curve} II.
 J. Reine Angew. Math. {\bf 233}
(1968), 89--100. 
\num{N}{M. Noether} {\"Uber die invariante Darstellung
algebraischer Funktionen}. Math. Ann. {\bf 17} (1880), 263--284.
\num{P}{K. Petri}{\"Uber die invariante Darstellung algebraischer
Funktionen einer Ver\"anderlichen}. Math. Ann. {\bf 88} (1923),
242--289.
\num{P-R}{K. Paranjape, S. Ramanan} {On the canonical ring of a
curve}. Algebraic geometry and commutative algebra, Vol. II,
503--516, Kinokuniya, Tokyo (1988). 
\num{S-D}{B. Saint-Donat}{On Petri's analysis of the linear system of
quadrics through a canonical curve}. Math. Ann. {\bf 206} (1973),
157--175. 
\num{S1}{F.-O. Schreyer} {Syzygies of canonical curves and
special linear series}. Math. Ann. {\bf 275} (1986),  105--137.
\num{S2}{F.-O. Schreyer}
{Green's conjecture for general p-gonal curves of large genus}.
Algebraic curves and projective geometry (Trento, 1988), 254--260,
Lecture Notes in Math. {\bf 1389}, Springer, Berlin, 1989.

\num{S3}{F.-O. Schreyer} {A standard basis approach to
syzygies of canonical curves}. J. Reine Angew. Math. {\bf 421} (1991),
83--123. 
\num{T}{M. Teixidor I Bigas}{Green's conjecture for the
generic $r$-gonal curve of genus} $g\geq 3r-7$. Duke Math. J.
{\bf 111} (2002),  195--222.
\num{V1}{C. Voisin}{Courbes t\'etragonales et
cohomologie de Koszul}.  J. Reine Angew. Math. {\bf 387} (1988),
111--121.
\num{V2}{C. Voisin}{Green's generic syzygy conjecture for curves of
even genus lying on a
$K3$ surface}. J. Eur. Math. Soc. {\bf 4} (2002),  363--404. 
\num{V3}{C. Voisin}{Green's canonical syzygy conjecture for
generic curves of odd genus}. Preprint {\tt math.AG/0301359}
(2003).
\vskip1cm\def\pc#1{\eightrm#1\sixrm}
\hfill\vtop{\eightrm\hbox to 5cm{\hfill Arnaud {\pc
BEAUVILLE}\hfill}
\hbox to 5cm{\hfill Laboratoire J.-A. Dieudonn\'e\hfill}
\hbox to 5cm{\sixrm\hfill UMR 6621 du CNRS\hfill}
\hbox to 5cm{\hfill  {\pc UNIVERSIT\'E DE} {\pc NICE}\hfill}
\hbox to 5cm{\hfill  Parc Valrose\hfill}
\hbox to 5cm{\hfill F-06108 {\pc NICE} Cedex 2\hfill}
\hbox to 5cm{\hfill \eightpoint\tt beauville@math.unice.fr\hfill}}\end
\end